\documentclass{amsproc}%
\usepackage{graphicx}
\usepackage{amscd}
\usepackage{amsmath}
\usepackage{amsfonts}
\usepackage{amssymb}%
\theoremstyle{plain}

\numberwithin{equation}{section}

\begin{document}
\title[Jordan normal base ]{The Jordan normal base in lattices and nilpotent endomorphisms of finitely
generated semisimple modules}
\author{Jen\H{o} Szigeti}
\address{Institute of Mathematics, University of Miskolc, Miskolc, Hungary 3515}
\email{jeno.szigeti@uni-miskolc.hu}
\thanks{\noindent This work was supported by OTKA of Hungary No. T043034 and K61007}

\begin{abstract}
We formulate a lattice theoretical Jordan normal form theorem for certain
nilpotent lattice maps satisfying the so called JNB conditions. As an
application of the general results, we obtain a transparent Jordan normal base
of a nilpotent endomorphism in a finitely generated semisimple module.

\end{abstract}
\maketitle

\noindent1. PRELIMINARIES

\bigskip

\noindent Let $\varphi:V\longrightarrow V$ be a nilpotent linear map of the
finite dimensional vector space $V$ over the complex number field $\mathbb{C}%
$. Now $\varphi^{k}=0\neq\varphi^{k-1}$ for some integer $1\leq k\leq n=\dim
V$ and we have a Jordan normal ($\mathbb{C}$-linear) base $\{v_{i}^{(t)}%
\mid1\leq t\leq r,1\leq i\leq k_{t}\}$ of $V$ with respect to $\varphi$ (see
[3]). The following conditions hold:

\noindent$k=\max\{k_{t}\mid1\leq t\leq r\}$, $k_{1}+k_{2}+...+k_{r}=n$ and%
\[
\varphi(v_{i}^{(t)})=v_{i-1}^{(t)}\text{ },\text{ }\varphi(v_{1}^{(t)})=0
\]
for all $1\leq t\leq r$ and $2\leq i\leq k_{t}$.

\noindent The purpose of this paper is to show, how the above classical result
is capable of broad generalization in the context of lattices. The present
work is of the same flavour as the job we carried out in [2]. We consider a
complete lattice $L$ and a nilpotent $\vee$-homomorphism $\lambda
:L\longrightarrow L$. The Jordan normal base of $L$ with respect to $\lambda$
is defined in a natural way. Using the properties of the induced map
$W\longmapsto\varphi(W)$ on the subspace lattice Sub$(V)$, we formulate the so
called JNB conditions. The main result of the paper is a lattice theoretical
Jordan normal form theorem (for nilpotent $\vee$-homomorphisms satisfying the
JNB conditions). Our proof is based on the use of the ideas presented in [4]
and [5]. The application of the general lattice theoretical results gives a
transparent Jordan normal base of a nilpotent endomorphism in a finitely
generated semisimple module. It is a natural challenge to take one further
step in the generalization process toward non-nilpotent maps.

\newpage

\noindent2. THE JORDAN NORMAL BASE IN A LATTICE

\bigskip

\noindent Let $(L,\vee,\wedge,0,1)$\ be a complete lattice with $\leq_{L}$
($<_{L}$) being the induced (strict) partial order on $L$ and consider a
complete $\vee$-homomorphism $\lambda:L\longrightarrow L$. The \textit{image}
and the \textit{kernel} of $\lambda$ can be defined as follows:%
\[
w=\text{im}(\lambda)=\lambda(1)\text{ and }z=\ker(\lambda)=\vee\{x\in
L\mid\lambda(x)=0\}.
\]
Clearly, the $\vee$-property of $\lambda$ ensures that $\lambda(x)\leq
_{L}\lambda(1)=w$ for all $x\in L$ and that $\lambda(x)=0\Longleftrightarrow
x\leq_{L}z$. If $\lambda^{k}=0\neq\lambda^{k-1}$ for some integer $1\leq k$
then $\lambda$ is called \textit{nilpotent} of index $k$. If $\lambda$ is
nilpotent and $x\leq_{L}\lambda(x)$ for some $x\in L$, then we have%
\[
x\leq_{L}\lambda(x)\leq_{L}\lambda^{2}(x)\leq_{L}...\leq_{L}\lambda
^{k}(x)=0_{L},
\]
whence $x=0$ (as well as $\lambda(0)=0$) follows. A finite set $\{a_{i}%
^{(t)}\mid1\leq t\leq r,1\leq i\leq k_{t}\}$ of atoms in $L$ is called a
\textit{Jordan normal base of }$L$ \textit{with respect to} $\lambda$, if we
have%
\[
a_{1}^{(1)}<_{L}a_{1}^{(1)}\vee a_{1}^{(2)}<_{L}...<_{L}a_{1}^{(1)}\vee
a_{1}^{(2)}\vee...a_{1}^{(r)}=c,
\]%
\[
c\vee a_{2}^{(1)}\vee...\vee a_{k_{1}}^{(1)}\vee...\vee a_{2}^{(t)}\vee...\vee
a_{i-1}^{(t)}<_{L}c\vee a_{2}^{(1)}\vee...\vee a_{k_{1}}^{(1)}\vee...\vee
a_{2}^{(t)}\vee...\vee a_{i-1}^{(t)}\vee a_{i}^{(t)},
\]%
\[
a_{1}^{(1)}\vee a_{2}^{(1)}\vee...\vee a_{k_{1}}^{(1)}\vee...\vee a_{1}%
^{(r)}\vee a_{2}^{(r)}\vee...\vee a_{k_{r}}^{(r)}=1_{L}%
\]
and%
\[
\lambda(a_{i}^{(t)})=a_{i-1}^{(t)}\text{ },\text{ }\lambda(a_{1}^{(t)})=0
\]
for all $1\leq t\leq r$ and $2\leq i\leq k_{t}$ ($a_{2}^{(t)}\vee...\vee
a_{i-1}^{(t)}$ is the empty join for $i=2$).

\bigskip

\noindent\textbf{2.1.Proposition.} \textit{If we have a Jordan normal base in
the complete lattice }$(L,\vee,\wedge,0,1)$\textit{\ with respect to the
complete }$\vee$\textit{-homomorphism }$\lambda:L\longrightarrow L$\textit{,
then }$\lambda$\textit{ is nilpotent.}

\bigskip

\noindent\textbf{Proof.} Let $\{a_{i}^{(t)}\mid1\leq t\leq r,1\leq i\leq
k_{t}\}$ be a Jordan normal base with respect to $\lambda$ and take
$k=\max\{k_{t}\mid1\leq t\leq r\}$. We claim, that $\lambda^{k}=0\neq
\lambda^{k-1}$. First we prove, that $\lambda^{k}(a_{i}^{(t)})=0$ for all
$1\leq t\leq r$, $1\leq i\leq k_{t}$. Clearly, we have%
\[
\lambda^{i}(a_{i}^{(t)})=\lambda^{i-1}(\lambda(a_{i}^{(t)}))=\lambda
^{i-1}(a_{i-1}^{(t)})=...=\lambda(a_{1}^{(t)})=0,
\]
whence $\lambda^{k}(a_{i}^{(t)})=\lambda^{k-i}(\lambda^{i}(a_{i}%
^{(t)}))=\lambda^{k-i}(0)=0$ can be deduced (note that $\lambda(0)=0$ is a
consequence of $0\leq_{L}a_{1}^{(1)}$). Since $\lambda^{k}$ is also a complete
$\vee$-homomorphism, we have%
\[
\lambda^{k}(1)=\vee\{\lambda^{k}(a_{i}^{(t)})\mid1\leq t\leq r,1\leq i\leq
k_{t}\}=0.
\]
Thus $\lambda^{k}=0$ and $\lambda^{k_{s}-1}(a_{k_{s}}^{(s)})=a_{1}^{(s)}\neq
0$, where $k=k_{s}$.$\square$

\bigskip

\noindent The following JNB conditions play a fundamental role in our
development, we formulate them for a pair $(L,\lambda)$, where $(L,\vee
,\wedge,0,1)$\ is a complete lattice and $\lambda:L\longrightarrow L$ is a
complete $\vee$-homomorphism.

\bigskip

\noindent JNB1: $L$ is of finite height and every element of $L$ is a finite
join of atoms. The height of $L$ is defined as%
\[
\text{height}(L)=\max\{m\in\mathbb{N}\mid x_{0}<_{L}x_{1}<_{L}...<_{L}%
x_{m}\text{ for some }x_{0},x_{1},...,x_{m}\text{ in }L\}.
\]

\bigskip

\noindent JNB2: For any choice of the elements $x\leq_{L}y$ with
$\lambda(x)=\lambda(y)$ we can find an element $u\in L$ such that $y=x\vee u$
and $\lambda(u)=0$.

\bigskip

\noindent JNB3: For any $x\in L$ the map $\lambda:[0,x]\longrightarrow
\lbrack0,\lambda(x)]$ is surjective.

\bigskip

\noindent If $x\leq_{L}y$ in $L$, then JNB1 implies that $y=a_{1}\vee...\vee
a_{m}$ for some atoms $a_{1},...,a_{m}$ in $L$. Now $y=x\vee a_{1}\vee...\vee
a_{m}$ and take $y=x\vee b_{1}\vee...\vee b_{s}$, where $b_{1},...,b_{s}$ are
atoms in $L$ with $s$ being the smallest possible. Then
\[
x<_{L}x\vee b_{1}<_{L}...<_{L}x\vee b_{1}\vee...\vee b_{s}=y
\]
is strictly ascending and $b_{1}\vee...\vee b_{s}$ is an irredundant join.

\bigskip

\noindent\textbf{2.2.Lemma.} \textit{JNB1 and JNB3 imply that for any atom
}$a\in\lbrack0,\lambda(1)]$\textit{ there exists an atom }$a^{\ast}\in
L$\textit{ such that }$a=\lambda(a^{\ast})$\textit{.}

\bigskip

\noindent\textbf{Proof.} Since $\lambda:[0,1]\longrightarrow\lbrack
0,\lambda(1)]$ is surjective, we can find an element $x\in L$ with
$a=\lambda(x)$. Now we have $x=b_{1}\vee b_{2}\vee...\vee b_{m}$, where each
$b_{i}\in L$ is an atom. Clearly, $b_{i}\leq_{L}x$ implies that $0\leq
_{L}\lambda(b_{i})\leq_{L}\lambda(x)=a$. For an index $i$, we have either
$\lambda(b_{i})=0$ or $\lambda(b_{i})=a$. If $\lambda(b_{i})=0$ holds for each
$1\leq i\leq m$, then we obtain%
\[
a=\lambda(x)=\lambda(b_{1}\vee b_{2}\vee...\vee b_{m})=\lambda(b_{1}%
)\vee\lambda(b_{2})\vee...\vee\lambda(b_{m})=0,
\]
a contradiction. Thus $\lambda(b_{j})=a$ holds for some $j\in\{1,2,...,m\}$%
.$\square$

\bigskip

\noindent We are ready to formulate the following converse of Proposition 2.1.

\bigskip

\noindent\textbf{2.3.Theorem.} \textit{Let }$(L,\vee,\wedge,0,1)$\ \textit{be
a complete lattice and }$\lambda:L\longrightarrow L$\textit{ a complete }%
$\vee$\textit{-homomorphism satisfying the above JNB conditions.} \textit{If
}$\lambda$\textit{ is nilpotent, then there exists a Jordan normal base in
}$L$\textit{ with respect to }$\lambda$\textit{.}

\bigskip

\noindent\textbf{Proof.} We apply an induction on the height of $L$. If
height$(L)=1$, then $L=\{0,1\}$ is a two element lattice and $\lambda(1)=0$
ensures that $\{1\}$ is a one element Jordan normal base of $L$ with respect
to $\lambda$ (note that $\lambda(1)=1$ would contradict to the nilpotent
property of $\lambda$). Let $n\geq2$ be an integer and assume that our theorem
holds for all lattices of height less or equal than $n-1$. Now consider the
situation described in the theorem with height$(L)=n$. Clearly, the interval
$L^{\ast}=[0,w]$ (with $w=\lambda(1)$) is a complete sublattice of $L$ and the
restriction $\lambda^{\ast}=\lambda\mid_{\lbrack0,w]}$ of $\lambda$ is a
complete $\vee$-homomorphism: $\lambda^{\ast}:L^{\ast}\longrightarrow L^{\ast
}$. It is straightforward to verify that the pair $(L^{\ast},\lambda^{\ast})$
satisfies all of the JNB conditions. In view of $\lambda(1)\neq1$, we have
height$(L^{\ast})\leq n-1$. Thus we have a Jordan normal base $\{a_{i}%
^{(t)}\mid1\leq t\leq r,1\leq i\leq k_{t}\}$ in $L^{\ast}$ with respect to
$\lambda^{\ast}$. Since $L^{\ast}=[0,\lambda(1)]$, Lemma 2.2 ensures that we
can find atoms $a_{k_{t}+1}^{(t)}\in L$, such that $\lambda(a_{k_{t}+1}%
^{(t)})=a_{k_{t}}^{(t)}$ holds for each $1\leq t\leq r$. Take $c=a_{1}%
^{(1)}\vee a_{1}^{(2)}\vee...\vee a_{1}^{(r)}$, then%
\[
a_{1}^{(1)}<_{L}a_{1}^{(1)}\vee a_{1}^{(2)}<_{L}...<_{L}a_{1}^{(1)}\vee
a_{1}^{(2)}\vee...a_{1}^{(r)}=c
\]
is a strictly ascending chain and%
\[
\lambda(c)=\lambda(a_{1}^{(1)}\vee a_{1}^{(2)}\vee...\vee a_{1}^{(r)}%
)=\lambda(a_{1}^{(1)})\vee\lambda(a_{1}^{(2)})\vee...\vee\lambda(a_{1}%
^{(r)})=0,
\]
whence $c\leq_{L}z=\ker(\lambda)$ follows. Using JNB1, we can exhibit a
sequence $b_{1},b_{2},...,b_{s}$ of atoms in $L$ such that%
\[
c<_{L}c\vee b_{1}<_{L}c\vee b_{1}\vee b_{2}<_{L}...<_{L}c\vee b_{1}\vee
b_{2}\vee...\vee b_{s}=z
\]
is a strictly ascending chain.

\noindent We claim, that $\{a_{i}^{(t)}\mid1\leq t\leq r,1\leq i\leq
k_{t}+1\}\cup\{b_{1},b_{2},...,b_{s}\}$ is a Jordan normal base in $L$ with
respect to $\lambda$. Now we have%
\[
z=a_{1}^{(1)}\vee a_{1}^{(2)}\vee...a_{1}^{(r)}\vee b_{1}\vee b_{2}\vee...\vee
b_{s},
\]
$b_{i}\leq_{L}z$, $1\leq i\leq s$ implies that $\lambda(b_{1})=\lambda
(b_{2})=...=\lambda(b_{s})=0$ and%
\[
z\vee a_{2}^{(1)}\vee...\vee a_{k_{1}+1}^{(1)}\vee...\vee a_{2}^{(t)}%
\vee...\vee a_{i-1}^{(t)}<_{L}z\vee a_{2}^{(1)}\vee...\vee a_{k_{1}+1}%
^{(1)}\vee...\vee a_{2}^{(t)}\vee...\vee a_{i-1}^{(t)}\vee a_{i}^{(t)}%
\]
holds for $1\leq t\leq r$ and $2\leq i\leq k_{t}+1$. Indeed,%
\[
z\vee a_{2}^{(1)}\vee...\vee a_{k_{1}+1}^{(1)}\vee...\vee a_{2}^{(t)}%
\vee...\vee a_{i-1}^{(t)}=z\vee a_{2}^{(1)}\vee...\vee a_{k_{1}+1}^{(1)}%
\vee...\vee a_{2}^{(t)}\vee...\vee a_{i-1}^{(t)}\vee a_{i}^{(t)},
\]
would imply that%
\[
a_{i}^{(t)}\leq_{L}z\vee a_{2}^{(1)}\vee...\vee a_{k_{1}+1}^{(1)}\vee...\vee
a_{2}^{(t)}\vee...\vee a_{i-1}^{(t)},
\]
whence first%
\[
a_{i-1}^{(t)}=\lambda(a_{i}^{(t)})\leq_{L}\lambda(z)\vee\lambda(a_{2}%
^{(1)})\vee...\vee\lambda(a_{k_{1}+1}^{(1)})\vee...\vee\lambda(a_{2}%
^{(t)})\vee...\vee\lambda(a_{i-1}^{(t)})=
\]%
\[
=a_{1}^{(1)}\vee...\vee a_{k_{1}}^{(1)}\vee...\vee a_{1}^{(t)}\vee...\vee
a_{i-2}^{(t)}%
\]
and then%
\[
c\vee a_{1}^{(1)}\vee...\vee a_{k_{1}}^{(1)}\vee...\vee a_{1}^{(t)}\vee...\vee
a_{i-2}^{(t)}=c\vee a_{1}^{(1)}\vee...\vee a_{k_{1}}^{(1)}\vee...\vee
a_{1}^{(t)}\vee...\vee a_{i-2}^{(t)}\vee a_{i-1}^{(t)}%
\]
can be derived, in contradiction with one of the properties of the Jordan
normal base $\{a_{i}^{(t)}\mid1\leq t\leq r,1\leq i\leq k_{t}\}$.

\noindent The use of JNB2 and%
\[
\lambda(a_{2}^{(1)}\vee...\vee a_{k_{1}+1}^{(1)}\vee...\vee a_{2}^{(r)}%
\vee...\vee a_{k_{r}+1}^{(r)})=
\]%
\[
=\lambda(a_{2}^{(1)})\vee...\vee\lambda(a_{k_{1}+1}^{(1)})\vee...\vee
\lambda(a_{2}^{(r)})\vee...\vee\lambda(a_{k_{r}+1}^{(r)})=
\]%
\[
=a_{1}^{(1)}\vee a_{2}^{(1)}\vee...\vee a_{k_{1}}^{(1)}\vee...\vee a_{1}%
^{(r)}\vee a_{2}^{(r)}\vee...\vee a_{k_{r}}^{(r)}=\lambda(1_{L})
\]
gives the existence of an element $u\in L$ such that%
\[
(a_{2}^{(1)}\vee...\vee a_{k_{1}+1}^{(1)}\vee...\vee a_{2}^{(r)}\vee...\vee
a_{k_{r}+1}^{(r)})\vee u=1_{L}%
\]
and $\lambda(u)=0$. Since $u\leq_{L}z$, we have%
\[
z\vee a_{2}^{(1)}\vee...\vee a_{k_{1}+1}^{(1)}\vee...\vee a_{2}^{(r)}%
\vee...\vee a_{k_{r}+1}^{(r)}=1_{L}%
\]
proving that $\{a_{i}^{(t)}\mid1\leq t\leq r,1\leq i\leq k_{t}+1\}\cup
\{b_{1},b_{2},...,b_{s}\}$ is a Jordan normal base in $L$ with respect to
$\lambda$.$\square$

\bigskip

\noindent The following is a trivial statement about the join of atoms.

\bigskip

\noindent\textbf{2.4.Proposition.} \textit{Let }$(L,\vee,\wedge,0)$%
\textit{\ be a lattice with the following properties}.

\noindent(i)\textit{ The cover relation }$x\vartriangleleft x\vee a$\textit{
holds for all }$x\in L$\textit{ and for all atoms }$a\in L$\textit{ with
}$x\neq x\vee a$\textit{ (atomic cover property, it is satisfied in upper
semimodular lattices).}

\noindent(ii)\textit{ Any two composition series of the form }$0=x_{0}%
\vartriangleleft x_{1}\vartriangleleft...\vartriangleleft x_{n}$\textit{ and}

\noindent$0=y_{0}\vartriangleleft y_{1}\vartriangleleft...\vartriangleleft
y_{m}$\textit{ with }$x_{n}=y_{m}$\textit{ are of the same length: }%
$n=m$\textit{.}

\noindent\textit{If }$a_{1},a_{2},...,a_{n}\in L$\textit{ are atoms such that}%
\[
a_{\pi(1)}<_{L}a_{\pi(1)}\vee a_{\pi(2)}<_{L}...<_{L}a_{\pi(1)}\vee a_{\pi
(2)}\vee...\vee a_{\pi(n)}%
\]
\textit{is a strictly ascending chain for some permutation }$\pi$\textit{ of
}$\{1,2,...,n\}$\textit{, then}

\noindent$a_{1}\vee a_{2}\vee...\vee a_{n}$\textit{ is an irredundant join (or
equivalently}%
\[
a_{1}\vee a_{2}\vee...\vee a_{n}=a_{1}\oplus a_{2}\oplus...\oplus a_{n}%
\]
\textit{is a direct sum).}

\newpage

\noindent3. NILPOTENT\ ENDOMORPHISMS\ OF\ FINITELY GENERATED

\noindent\ \ \ \ SEMISIMPLE MODULES

\bigskip

\noindent Let $R$ be a ring and $\varphi:M\longrightarrow M$ an $R$%
-endomorphism of the unitary left $R$-module $_{R}M$. For a submodule $N\leq
M$ let $\varphi(N)=\{\varphi(x)\mid x\in N\}$ denote the $\varphi$-image of
$N$. Clearly, $\varphi(N)\leq M$ is an $R$-submodule and for any family
$N_{\gamma}$, $\gamma\in\Gamma$ of $R$-submodules in $M$ we have%
\[
\varphi\left(  \underset{\gamma\in\Gamma}{\sum}N_{\gamma}\right)
=\underset{\gamma\in\Gamma}{\sum}\varphi(N_{\gamma}).
\]
If $X,Y\leq M$ are $R$-submodules with $X\subseteq Y$ and $\varphi
(X)=\varphi(Y)$, then%
\[
X+(Y\cap\ker(\varphi))=Y.
\]
Take $y\in Y$, then $\varphi(y)=\varphi(x)$ for some $x\in X$ and $y=x+(y-x)$
with $\varphi(y-x)=\varphi(y)-\varphi(x)=0$. Thus $y-x\in Y\cap\ker(\varphi)$
and our claim is proved. We note that $\varphi(Y\cap\ker(\varphi))=\{0\}$
holds for the $R$-submodule $Y\cap\ker(\varphi)\leq M$.

\noindent If $X\leq M$ is an $R$-submodule with $X\subseteq\varphi(N)$, then
consider%
\[
\varphi^{-1}(X)=\{u\in M\mid\varphi(u)\in X\}
\]
an $R$-submodule of $M$ and take $Z=N\cap\varphi^{-1}(X)\subseteq N$. For an
element $x\in X$, the containment $X\subseteq\varphi(N)$ ensures that we can
find an element $z\in N$ such that $x=\varphi(z)$. Now $x\in X$ implies that
$z\in N\cap\varphi^{-1}(X)$, whence $\varphi(Z)=X$ follows.

\noindent In view of the above observations, we have JNB2 and JNB3 for the
pair $(L,\lambda)$, where%
\[
L=\text{Sub}(_{R}M)
\]
is the complete lattice of the $R$-submodules of $M$ with respect to the
containment relation $\subseteq$\ and%
\[
\lambda(N)=\varphi(N)
\]
is an $L\longrightarrow L$ complete $\vee$-homomorphism (we note that
$\underset{\gamma\in\Gamma}{\vee}N_{\gamma}=\underset{\gamma\in\Gamma}{\sum
}N_{\gamma}$).

\noindent The following conditions are equivalent for the modular lattice
$L=$Sub$(_{R}M)$:

\begin{enumerate}
\item $L=$Sub$(_{R}M)$ satisfies JNB1.

\item $_{R}M$ is a semisimple left $R$-module of finite height.

\item $_{R}M$ is a finitely generated semisimple left $R$-module.
\end{enumerate}

\bigskip

\noindent Now we are ready to present the Jordan normal base theorem for
nilpotent $R$-module endomorphisms.

\bigskip

\noindent\textbf{3.1.Theorem.}\textit{ Let }$\varphi:M\longrightarrow
M$\textit{ be a nilpotent }$R$\textit{-endomorphism of the finitely generated
semisimple left }$R$\textit{-module }$_{R}M$\textit{. Then there exists a
subset}

\noindent$\{x_{i}^{(t)}\mid1\leq t\leq r,1\leq i\leq k_{t}\}$\textit{ in }%
$M$\textit{ such that each }$R$\textit{-submodule }$Rx_{i}^{(t)}\leq
M$\textit{ is simple,}%
\[
\underset{1\leq t\leq r,1\leq i\leq k_{t}}{\oplus}Rx_{i}^{(t)}=M
\]
\textit{is a direct sum and }$\varphi(x_{i}^{(t)})=x_{i-1}^{(t)}$\textit{ }%
$,$\textit{ }$\varphi(x_{1}^{(t)})=0$\textit{ for all }$1\leq t\leq r$\textit{
and }$2\leq i\leq k_{t}$\textit{.}

\bigskip

\noindent\textbf{Proof.} Consider the pair $(L,\lambda)$, where $L=$%
Sub$(_{R}M)$ and $\lambda(N)=\varphi(N)$ for an $R$-submodule $N\leq M$ of
$M$. As we have already seen, $(L,\lambda)$ satisfies the conditions JNB1,
JNB2 and JNB3. We note that any atomistic (upper semi-) modular lattice of
finite height (such as $L$) satisfies the conditions (i) and (ii) of
Proposition 2.4. Since $\varphi^{k}=0\neq\varphi^{k-1}$ implies that
$\lambda^{k}=0\neq\lambda^{k-1}$, the complete $\vee$-homomorphism $\lambda$
is nilpotent. Thus we can apply Theorem 2.3 to obtain a Jordan normal base
$\{V_{i}^{(t)}\mid1\leq t\leq r,1\leq i\leq k_{t}\}$ of $L$ with respect to
$\lambda$. Clearly, each $V_{i}^{(t)}\leq M$ is a simple $R$-submodule of $M$
and Proposition 2.4 guarantees that%
\[
\underset{1\leq t\leq r,1\leq i\leq k_{t}}{\oplus}V_{i}^{(t)}=M
\]
is a direct sum. Since%
\[
\varphi(V_{i}^{(t)})=V_{i-1}^{(t)}\text{ },\text{ }\varphi(V_{1}^{(t)})=0
\]
for all $1\leq t\leq r$ and $2\leq i\leq k_{t}$, starting with an arbitrary
$0\neq x_{1}^{(t)}\in V_{1}^{(t)}$ element, we can take the consecutive
$\varphi$-preimages to obtain a sequence $x_{i}^{(t)}\in V_{i}^{(t)}$, $1\leq
i\leq k_{t}$\ such that%
\[
\varphi(x_{i}^{(t)})=x_{i-1}^{(t)}\text{ (and }\varphi(x_{1}^{(t)})=0\text{)}%
\]
for all $1\leq t\leq r$ and $2\leq i\leq k_{t}$. To conclude the proof, it is
enough to note that $x_{i}^{(t)}\neq0$ in the simple $R$-module $V_{i}^{(t)}$,
whence $Rx_{i}^{(t)}=V_{i}^{(t)}$ follows.$\square$

\bigskip

\bigskip

\noindent REFERENCES

\bigskip

\begin{enumerate}
\item Anderson, F.W. and Fuller, K.R.: \textit{Rings and Categories of
Modules}, GTM Vol. 13, Springer Verlag, New York, 1974.

\item K\"{o}rtesi, P. and Szigeti, J.:\textit{ A general approach to the
Fitting lemma}, Mathematika (London) Vol. 52 (2005), 155-160.

\item Prasolov, V.V.: \textit{Problems and Theorems in Linear Algebra}, Vol.
134 of Translation of Mathematical Monographs, American Mathematical Society,
Providence, Rhode Island, 1994.

\item Valiaho, H.: \textit{An elementary approach to the Jordan normal form of
a matrix}, Amer. Math. Monthly Vol. 93 (1986), 711-714.

\item Wildon, M.: \textit{Jordan normal form and the structure theorem for
abelian groups}, (manuscript, January 2005).
\end{enumerate}

\end{document}